\documentclass[11pt,reqno]{amsart}
\usepackage[margin=1in]{geometry}
\usepackage{algorithm,algpseudocode}
\usepackage[fleqn,tbtags]{mathtools}
\usepackage{mathtools}
\usepackage{pifont,enumerate}
\usepackage{verbatim}
\usepackage{tikz}

\usepackage{filecontents}
\usepackage{placeins}
\usepackage{mlmodern}
\usepackage{extarrows}
\usepackage{graphicx}
\usepackage{calc}
\usepackage{adjustbox}
\usepackage{amssymb}

\usepackage{subcaption}
\usepackage{graphicx}
\usepackage[american, siunitx]{circuitikz}
\usetikzlibrary{graphs}
\usetikzlibrary{quotes,arrows.meta}

\usepackage{mathtools}

\usepackage{hyperref,xcolor}
\hypersetup{
    colorlinks,
    linkcolor={red!50!black},
    citecolor={blue!50!black},
    urlcolor={blue!80!black}
}

\DeclarePairedDelimiterX\innerp[2]{\langle}{\rangle}{#1,#2}

\newcommand\restr[2]{{
  \left.\kern-\nulldelimiterspace 
  #1 
  \vphantom{\big|} 
  \right|_{#2} 
  }}

\theoremstyle{definition}
\newtheorem{theorem}{Theorem}[section]

\newtheorem{example}[theorem]{Example}
\numberwithin{equation}{section}

\newcommand*{\method}[1]{#1}
\newcommand*{\overeqU}[1]{\ensuremath{\mathrel{\overset{\method{#1}}{=}}}}
\newcommand*{\overeq}[1]{\mathrel{\overset{\method{#1}}{\resizebox{\widthof{\kern1.25pt\overeqU{\method{#1}}}}{\heightof{$=$}}{$=$}}}}

\DeclareMathOperator{\vol}{Vol}

\newcommand\cuboid[7]{
\draw[black,fill= #7] (#1+#4/2,#2+#5/2,#3+#6/2) -- ++(-#4,0,0) -- ++(0,-#5,0) -- ++(#4,0,0) -- cycle;
\draw[black,fill=#7] (#1+#4/2,#2+#5/2,#3+#6/2) -- ++(0,0,-#6) -- ++(0,-#5,0) -- ++(0,0,#6) -- cycle;
\draw[black,fill=#7] (#1+#4/2,#2+#5/2,#3+#6/2) -- ++(-#4,0,0) -- ++(0,0,-#6) -- ++(#4,0,0) -- cycle;
}

\begin{document}

\title{Geometric Programming for 3D Circuits}

\author[R.~Wang]{Rongbiao~Wang}
\address{Committee on Computational and Applied Mathematics, University of Chicago, Chicago, IL 60637}
\email{rbwang@uchicago.edu}
\author[L.-H.~Lim]{Lek-Heng~Lim}
\address{Computational and Applied Mathematics Initiative, Department of Statistics,
University of Chicago, Chicago, IL 60637-1514}
\email{lekheng@uchicago.edu}

\begin{abstract}
    With the soaring demand for high-performing integrated circuits, 3D integrated circuits (ICs) have emerged as a promising alternative to traditional planar structures. Unlike existing 3D ICs that stack 2D layers, a full 3D IC features cubic circuit elements unrestricted by layers, offering greater design freedom. Design problems such as floorplanning, transistor sizing, and interconnect sizing are highly complex due to the 3D nature of the circuits and unavoidably require systematic approaches. We introduce geometric programming to solve these design optimization problems systematically and efficiently. 
\end{abstract}

\maketitle

\section{Introduction}

In recent years, improving the energy efficiency of 2D integrated circuits (ICs) has become increasingly difficult due to physical, lithographic, manufacturing, and energy constraints \cite{power_benefit,Moore}. As demand for higher-performing chips continues to increase, 3D integration has become a promising candidate to tackle these limitations.

Current 3D integration succeeds by stacking layers of semiconductor components on top of each other to create smaller and more effective electronic devices and to integrate heterogeneous functionalities. However, the limited design freedom in these layer-based approaches reduces their adaptability to a broader range of applications \cite{guha2022future}. A full 3D IC eliminates layer-based constraints by employing truly three-dimensional components without being confined to layers: The circuit elements are no longer constrained by height and can be oriented and placed in a three-dimensional structure. In both electronic and photonic integrated circuits, full 3D integration allows for further space minimization and higher integrability. Moreover, full 3D ICs can realize any given circuit topology without crossing of wires, while 2D ICs can only do so for a small subclass of circuits. In Section~\ref{sec:IC}, we provide an overview of full 3D integration, laying the foundation for formulating the design problems.

The design problems aim to optimize the trade-offs among volume, energy efficiency, delay, and other key metrics. These problems--such as floorplanning, gate sizing, and interconnect sizing--must be solved across multiple abstraction levels, even in traditional 2D integrations. Given that modern circuits contain billions of transistors, manual computation is infeasible, necessitating a systematic approach. 

In Section~\ref{sec:design}, we show that the geometric programming approach of Boyd, Kim, Patil, and Horowitz \cite{digital_gp} naturally extends to the design of 3D ICs. Nevertheless, there is an important distinction: We argue that while the design of 2D chips can be done without such a sophisticated approach, the design of 3D chips are of such a level of complexity that makes computer-aided approaches like geometric programming all but inevitable. In other words, we think that the approach of Boyd et al.\ will likely have a greater impact on 3D chip design. We show here that it essentially applies to 3D chip with only minor modifications. To our knowledge, this is the first work to bring the GP-based approach to the context of full 3D integration.

\section{Full 3D integration}\label{sec:IC}

We start by outlining two key insights from current technologies that motivate the need for full 3D integration.

The first insight comes from the current 3D ICs, most notably TSV-based ICs and monolithic 3D ICs. Known as chip stacking, TSV-based ICs vertically connect multiple semiconductor device wafers using through-silicon vias (TSVs) \cite{stacking1,stacking2}. In contrast, monolithic 3D ICs directly connect functional device layers through inter-tier vias (MIVs) without placing them on separate wafers \cite{monolithic1,monolithic2}. 

The lesson is that dividing an IC into layers is suboptimal. When each layer resides on a separate wafer, as in chip stacking, TSVs are required to connect the wafers. These micrometer-scaled TSVs restrict space optimization due to their large size compared to nanometer-scaled transistors \cite{Monolithic-vs-TSV}. Additionally, they demand complex wafer fabrication and bonding processes, which significantly constrain integrability and adaptability \cite{monolithic_solution}. Monolithic integration enhances functionality in heterogeneous integration by eliminating the need for separate wafers \cite{guha2022future}. However, even in monolithic 3D ICs, confining circuit elements to layers still limits design flexibility. This layer-based approach necessitates custom engineering for each specific device, reducing the adaptability of these technologies. These limitations motivate full 3D integration, which removes the layer constraint entirely.

\begin{figure}
    \centering
    \begin{subfigure}{0.45\textwidth}
    \centering
    \includegraphics[trim={6em 8ex 1em 8ex},clip,width=0.95\textwidth]{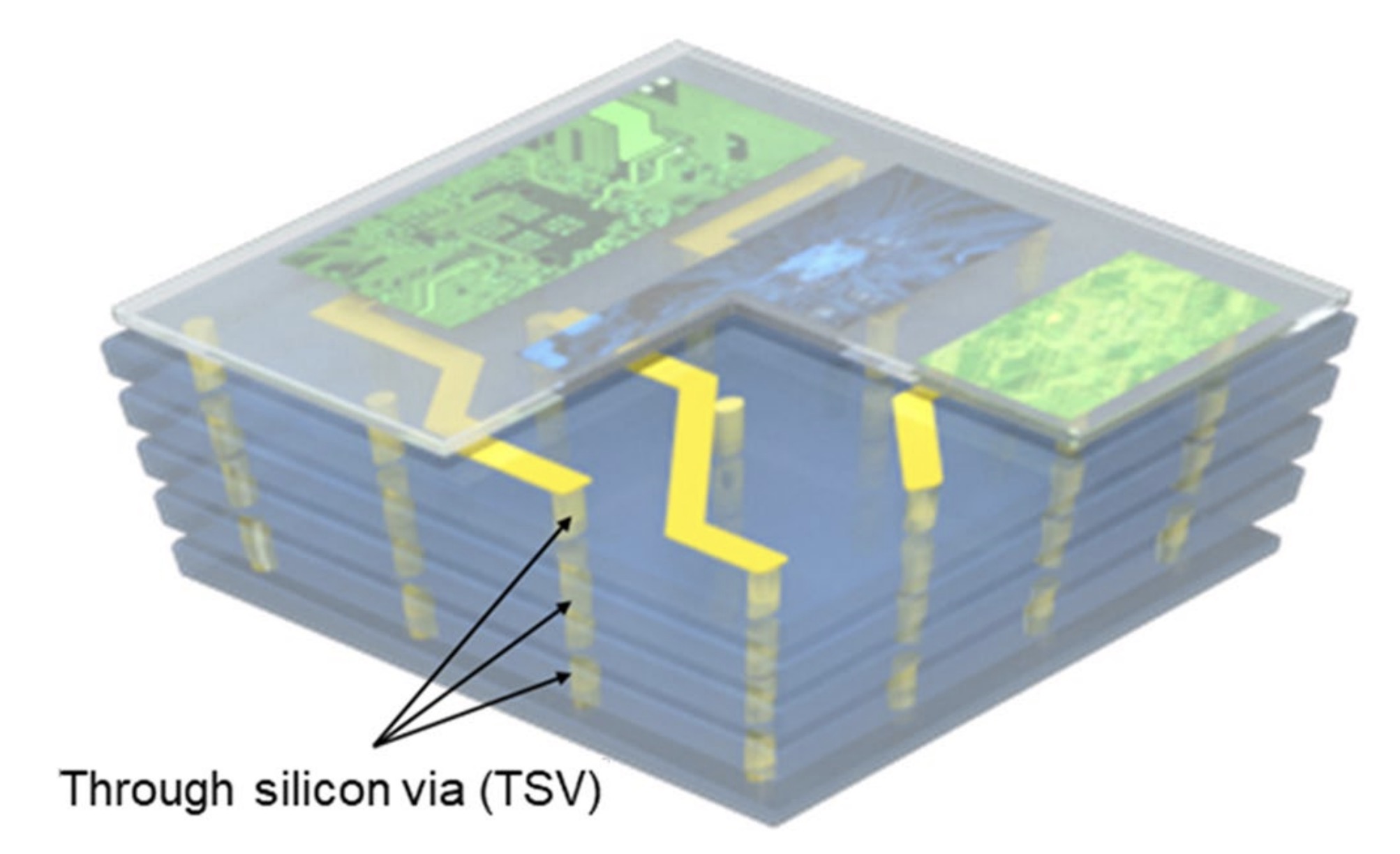}
    \caption{A TSV-based 3D IC \cite{TSV_image}.}
      \end{subfigure}
    \begin{subfigure}{0.54\textwidth}
        \centering
      \includegraphics[trim={3em 6.5ex 7em 4ex},clip,width=0.9\textwidth]{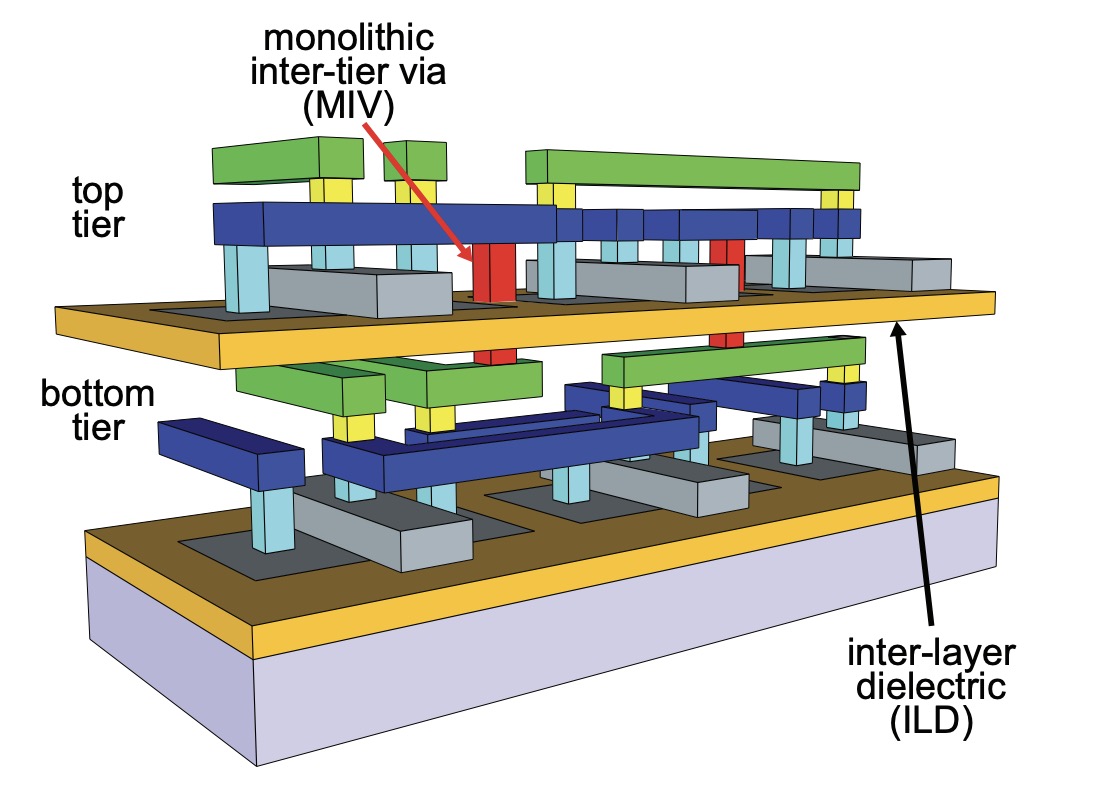}
      \caption{A monolithic 3D IC \cite{monolithic_image}.}
    \end{subfigure}
    \caption{Chip stacking and monolithic 3D IC}
\end{figure}

The second insight comes from examining the inherently 3D physical structure of modern transistors. Since the early 2000s, 3D transistors such as FinFETs \cite{FinFET_original,FinFET} and gate-all-around FETs (GAAFETs) \cite{multi-gate} have largely replaced traditional metal-oxide-semiconductor field-effect transistors (MOSFETs) in ICs. Additionally, tunneling FETs (TFETs), another alternative to MOSFETs, have recently demonstrated significant potential for leveraging 3D geometries in their design \cite{tunnelling-transistors}. These advancements suggest that transistors are better modeled as cubic objects rather than flat components constrained by height in 3D ICs.

\begin{figure}[htb]
  \includegraphics[trim={1em 3ex .5em 2ex},clip,width=0.6\textwidth]{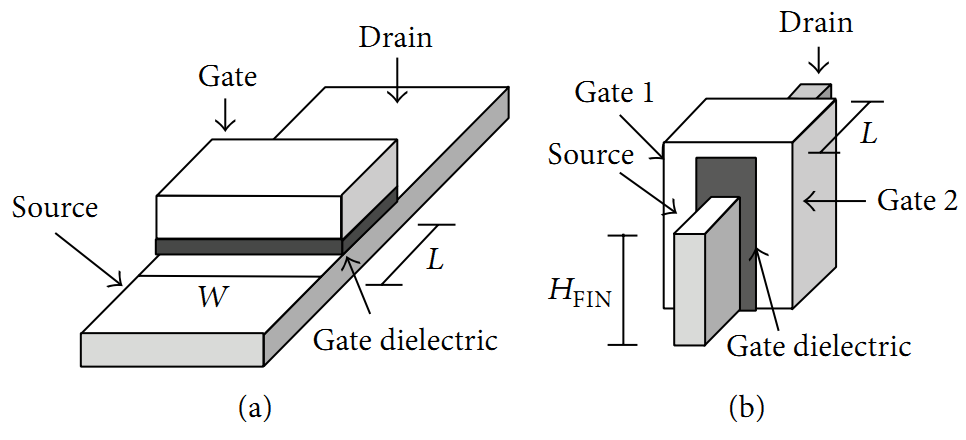}
  \caption{Planar transistor vs. FinFET \cite{FinVMosPic}.} \label{fig:finfet}
\end{figure}

\begin{figure}[tbh]
    \includegraphics[trim={12em 16ex 13em 14ex},clip,width=0.6\textwidth]{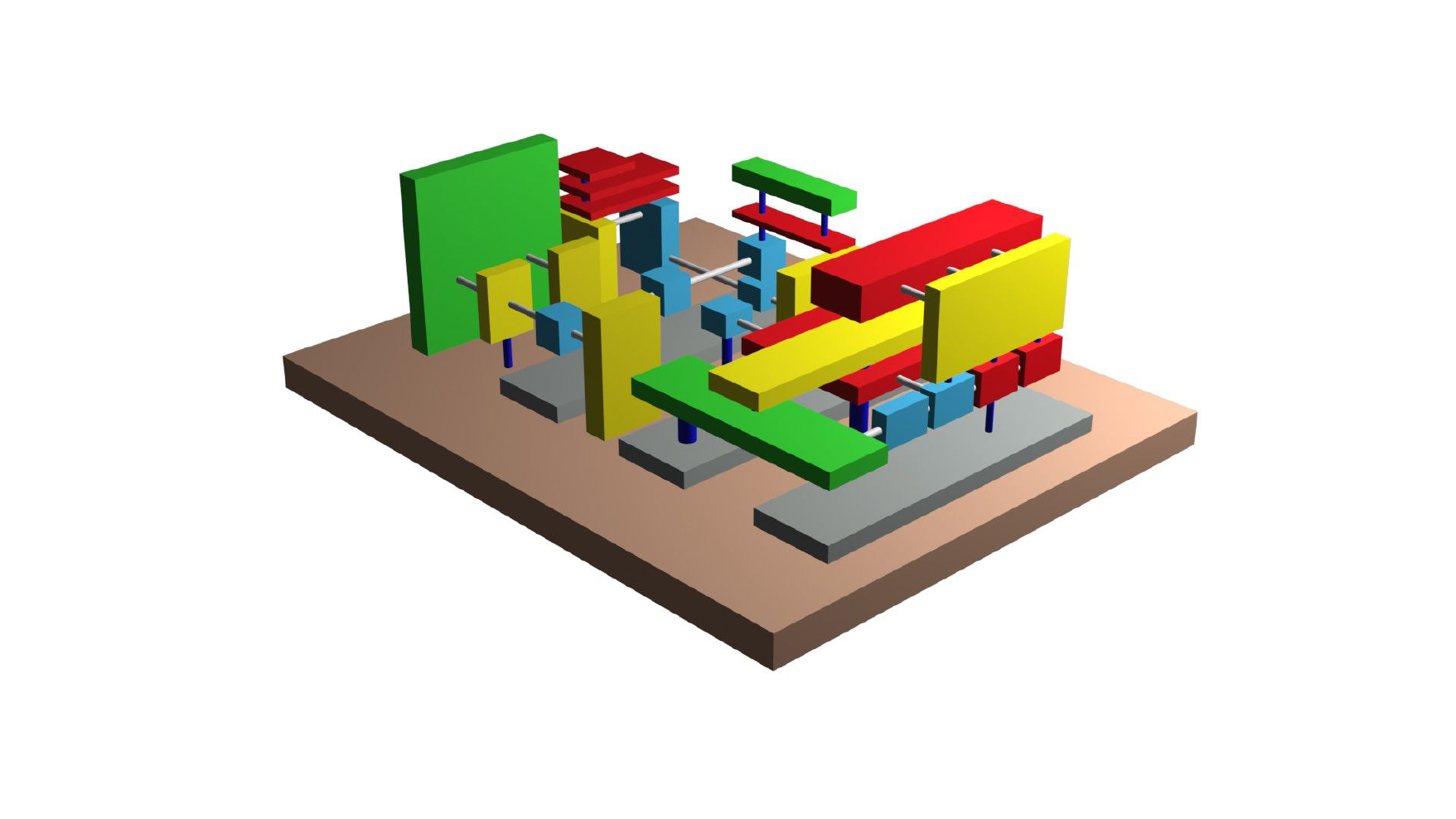}
    \caption{A full 3D IC.} \label{fig:full_3d}
\end{figure}

Building on these insights, we propose a full 3D integrated circuit with the following characteristics:
\begin{enumerate}
  \item Circuit elements are not confined to specific layers;
  \item Each circuit element is a cubical object with comparable $x$, $y$, and $z$ dimensions;
  \item Interconnects directly connect individual elements.
\end{enumerate}
The term ``circuit elements" here is intentionally broad, encompassing transistors, gates, blocks, or other structures, depending on the level of integration and design. Figure~\ref{fig:full_3d} illustrates an example of a full 3D IC.

The full 3D integration can be applied to both 3D electronic circuits, which transmit information via electrons, and 3D photonic integrated circuits (PICs) \cite{photonic_1,photonic_2}, which use photons. This flexibility allows for further exploration of PICs, which can overcome the heat dissipation challenges of electronic circuits. Fabrication of full 3D integrations poses significant technological barriers as their complex geometries introduce additional challenges. 3D printing \cite{3Dprinting2,3Dprinting3} provides a promising solution due to its capacity to handle intricate structures. Recent research has suggested it as a promising direction for future development in 3D IC fabrications \cite{3Dprinting1}.

\section{Design Problems}\label{sec:design}

Design optimization takes place once the relative positions of the circuit elements and interconnects are established. Determining these positions corresponds to finding a 3D graph embedding of the circuit topology, a well-studied topic with established algorithms \cite{3d_draw1,3d_draw6,3d_draw3,3d_draw5,3d_draw6,3d_draw7}. In fact, these algorithms reveal a significant but often overlooked advantage of 3D integration: any circuit topology could be embedded in 3D without interconnect crossings. In contrast, a crossing-free 2D embedding is possible only for planar graphs. As shown in \cite{planar_limit}, the fraction of planar graphs among all graphs with $n$ vertices approaches zero as $n$ grows:
\[
  \lim_{n \to \infty}\frac{\#\{\text{planar graphs with $n$ vertices}\}}{\#\{\text{graphs with $n$ vertices}\}} = 0.
\]
Consequently, as circuit complexity increases, the likelihood of embedding a circuit topology in 2D without crossings becomes negligible. As a result, planar integrated circuits with nonplanar circuit topologies require interconnects to bridge over one another, leading to longer connection lengths and increased communication overhead. 

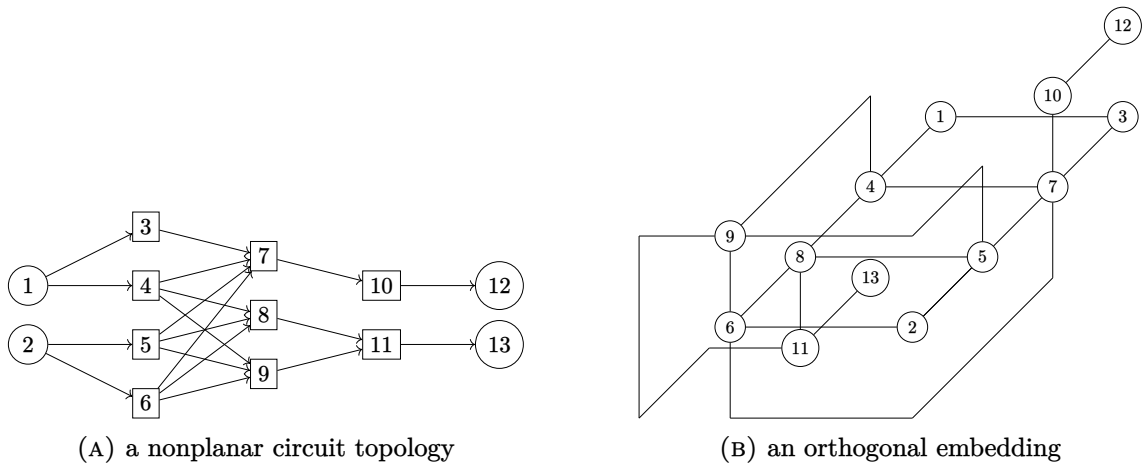
\begin{figure}
    \centering
    \begin{subfigure}{0.45\textwidth}
    \centering
    \resizebox{0.95\textwidth}{!}{
        \begin{tikzpicture}
\node [circle,draw] (1) at (0,0) {1};
\node [circle,draw] (2) at (0,-1) {2};
\node [draw] (3) at (2,1) {3};
\node [draw] (4) at (2,0) {4};
\node [draw] (5) at (2,-1) {5};
\node [draw] (6) at (2,-2) {6};
\node [draw] (7) at (4,0.5) {7};
\node [draw] (8) at (4,-.5) {8};
\node [draw] (9) at (4,-1.5) {9};
\node [draw] (10) at (6,0) {10};
\node [draw] (11) at (6,-1) {11};
\node [circle,draw] (12) at (8,0) {12};
\node [circle,draw] (13) at (8,-1) {13};   
\graph { (1) -> {(3), (4)};
(2) -> {(5), (6)};
(3) -> (7);
(5) -> {(7), (8), (9)} -> {(10), (11)} -> {(12), (13)};
(4) -> {(7), (8), (9)};
(6) -> {(7), (8), (9)};
};
\end{tikzpicture}
    }
      \caption{a nonplanar circuit topology}
      \label{fig:nonplanar}
      \end{subfigure}
    \begin{subfigure}{0.54\textwidth}
        \centering
      \resizebox{0.75\textwidth}{!}{\begin{tikzpicture}
    \node [circle,draw] (1) at (0,2,-4){$1$};
    \node [circle,draw] (2) at (4,2,8){$2$};
    \node [circle,draw] (3) at (4,2,-4){$3$};
    \node [circle,draw] (4) at (0,2,0){$4$};
    \node [circle,draw] (5) at (4,2,4){$5$};
    \node [circle,draw] (6) at (0,2,8){$6$};
    \node [circle,draw] (7) at (4,2,0){$7$};
    \node [circle,draw] (8) at (0,2,4){$8$};
    \node [circle,draw] (9) at (0,4,8){$9$};
    \node [circle,draw] (10) at (4,4,0){$10$};
    \node [circle,draw] (13) at (0,0,0){$13$};
    \node [circle,draw] (12) at (4,4,-4){$12$};
    \node [circle,draw] (11) at (0,0,4){$11$};
    \path {
        (1) edge (4) (4) edge (8) (8) edge (6) (6) edge (2) (2) edge (5) (5) edge (7) (4) edge (7) (8) edge (5) (2) edge (5) (7) edge (10) (10) edge (12) (6) edge (9) (1) edge (3) (3) edge (7) (11) edge (13) (8) edge (11)
        (6) edge ++(0,-2, 0) (0,0,8) edge ++(4,0,0) (4,0,8) edge ++(0,0,-8) (4,0,0) edge (7)
        (9) edge ++(4,0,0) (4,4,8) edge ++(0,0,-4) (4,4,4) edge (5)
        (9) edge ++(0,0,-8) (0,4,0) edge (4)
        (9) edge ++(-2,0,0) (-2,4,8) edge ++(0,-4,0) (-2,0,8) edge ++(0,0,-4) (-2,0,4) edge (11)
        };
    \end{tikzpicture}}
      \caption{an orthogonal embedding}
      \label{fig:ortho}
    \end{subfigure}
    \caption{A nonplanar graph and its orthogonal embedding. }
    \label{fig:graph}
\end{figure}

Optimizing the performance of an IC under various constraints is a complex task. It requires balancing trade-offs among volume, temperature, energy efficiency, and other factors. This challenge becomes even more intricate in full 3D integration: The cubic structure of circuit elements offers greater optimization potential, but also introduces a larger set of variables, making the optimization process significantly more complex. We propose using geometric programming (GP), a widely applied class of optimization problems in engineering \cite{gas,PowerControl,aircraft}, as a systematic approach to formulating and solving these problems.

\subsection{Geometric programming}\label{sec:GP}

Let $x = (x_1,\dots,x_n) \in \mathbb{R}^n_+$ denote the optimization variables for the rest of Section~\ref{sec:design}. A \emph{posynomial} is a real-valued function $f$ on $\mathbb{R}^n_+$ of the form
\[
  f(x) = \sum_{i=1}^m c_i x_1^{\alpha_{i1}}\dots x_n^{\alpha_{in}}
\]
with $c_i > 0$ and $\alpha_{ij} \in \mathbb{R}$, $i=1,\dots,m$, $j = 1,\dots,n$. When $m=1$, the function $f$ is called a \emph{monomial}. The posynomials and monomials here differ from polynomials and monomials in standard algebra by having positive coefficients and allowing fractional exponents instead of having real coefficients and nonnegative integer exponents. A \emph{geometric program} is a constrained optimization problem of the form 
\begin{equation}\label{eq:GP}
\begin{aligned}
    \text{minimize} \quad & f_0(x) \\
    \text {subject to} \quad & f_i(x) \leq 1, \quad i=1, \dots, p, \\
    & g_j(x)=1, \quad j=1, \dots, q,
\end{aligned}
\end{equation}
where $f_i$ are posynomials for $i = 0,\dots,p$ and $g_j$ are monomials for $j = 1,\dots,q$. 

Many optimization problems are extensions of GP, meaning they can be transformed into and solved as an equivalent GP. An important class of extensions is \emph{generalized geometric program} (GGP), which replaces the posynomials in \eqref{eq:GP} by \emph{generalized posynomials}. Generalized posynomials are functions formed from additions, multiplications, positive powers, and maximums of posynomials. In the following sections, we will formulate three design problems as GGPs, which can be efficiently solved as GPs without loss of generality.

A key advantage of GP modeling is the existence of efficient and readily available software packages. The standard approach to solve \eqref{eq:GP} is to convert it into the problem
\begin{equation}\label{eq:convex_GP}
  \begin{aligned}
      \text{minimize} \quad & \log f_0(e^y) \\
  \text { subject to} \quad & \log f_i(e^y) \leq 0, \quad i=1, \dots, p, \\
  & \log g_j(e^y)=0, \quad j=1, \dots, q,
  \end{aligned}
\end{equation}
where $y = \log x = (\log x_1,\dots, \log x_n)$ and $e^y = (e^{y_1},\dots,e^{y_n})$. Since the transformed functions $\log f_i(e^y)$ are convex for $i=0,\dots,p$ and $\log g_j(e^y)$ are affine for $j = 1, \dots, q$, problem \eqref{eq:convex_GP} is a convex optimization problem. Algorithms like the interior-point method efficiently solve large-scale convex problems. This approach is already implemented in many existing convex optimization software packages such as \texttt{GGPLAB} in \textsc{Matlab} and \texttt{CVXPY} in Python, allowing users to directly optimize a GGP. Compared to other numerical optimization methods, GP-based methods do not require algorithm parameter tuning and consistently obtain the global solution.

Another advantage of GP-based methods is their efficiency in finding parameters in objective and constraint functions from experimental data. A function $f(x)$ can be approximated by a monomial (resp. a posynomial) if $\log f(e^y)$ is nearly affine (resp. convex). When these conditions hold, fitting parameters from data reduces to solving a nonlinear least-squares problem using the Gauss-Newton algorithm. Interested readers may refer to \cite{tutorial} for a more detailed introduction to GP. Now we introduce three design problems in full 3D IC formulated as GPs. All codes in the following sections have been made available at:
\begin{quote}
    \url{https://github.com/thomasw15/GP_3D/tree/main}.
\end{quote}

\subsection{Floorplanning}\label{sec:floorplanning}
Among design challenges, thermal management is particularly critical in 3D ICs, addressed through innovations such as microfluidic cooling and heat spreaders \cite{heat_transfer,spreader}. Although minimizing the size of an IC is desirable, it must be balanced with effective heat dissipation. Floorplanning aims to achieve this balance within size constraints. Based on \cite{temperature}, we propose a temperature-aware GP formulation for floorplanning in full 3D ICs. Although the term ``floorplanning" is technically a misnomer in the context of 3D integration, where layers or floors are absent, we use it for consistency with established terminology.

Floorplanning assumes that the relative positions of the integrated circuit modules, including circuit elements and heat removal technologies, are predetermined. Suppose there are $n$ modules indexed by $1,\dots,n$, each with dimensions $(x_i, y_i, z_i)$, which are the optimization variables. Let $X$, $Y$, and $Z$ denote the dimensions of the smallest cube that encloses all modules and $T_i$ the temperature of module $i$. The objective is to minimize
\[\alpha XYZ + (1-\alpha)\sum_{i=1}^n T_i,\]
where the constant $0 \leq \alpha \leq 1$ balances size and temperature. If the IC is a PIC where heat is not a concern, we set $\alpha = 1$. 

If module $i$ is a circuit element, its temperature is 
\[T_i = P_i K_i^{-1}t_ia_i^{-1}\]
where $t_i$ is the thickness, $a_i$ is the area of its flat face, and $P_i$ and $K_i$ are given constants for power consumption and thermal conductivity, respectively. Both thickness and area depend on the orientation of the module in space as shown in Example~\ref{ex:floorplanning}. For a heat removal module $i$, we set $T_i = 0$.

The floorplanning is subject to the following constraints. First,
\[
  x_i \geq x_i^{\min}, \quad y_i \geq y_i^{\min}, \quad z_i \geq z_i^{\min}, \quad i = 1,\dots, n,
\] 
where $x_i^{\min}$, $y_i^{\min}$, and $z_i^{\min}$ are fixed minimal dimensions of the modules. Second, the total thickness is limited by
\[
  Z \leq Z^{\max}
\]
due to manufacturing constraints. Finally, the relative positions of the modules define how their dimensions contribute to the overall dimensions of the bounding cube. This is best illustrated through the following example.

\begin{example}\label{ex:floorplanning}
Consider four transistor modules arranged as shown in Figure~\ref{fig:floorplanning}. The dimensional constraint in the $X$-direction is 
\[\max\{x_1+x_4,x_2+x_4,x_3+x_4\} \leq X\]
since modules $1$, $2$, and $3$ are aligned adjacent to $4$ in the $X$-direction. Similarly, the constraints in the other two directions are
\[\max\{y_1,y_2,y_3,y_4\} \leq Y,\]
and
\[\max\{z_1+z_2+z_3,z_4\}\leq Z.\]

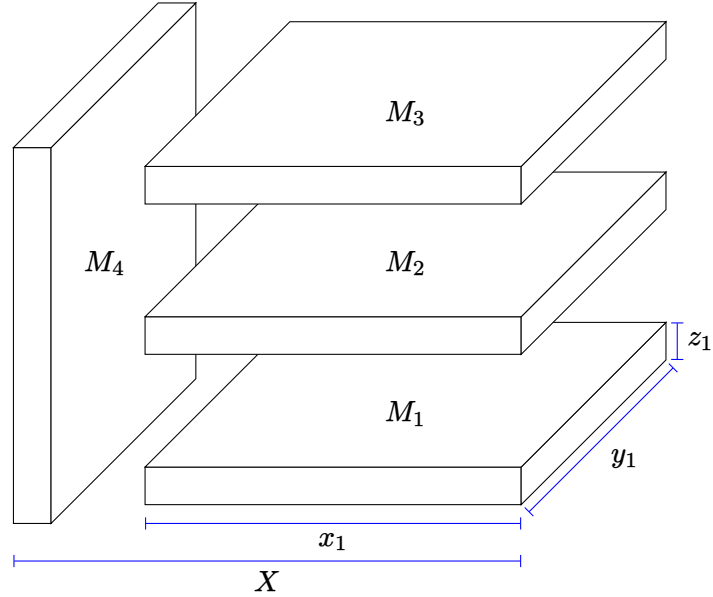
\begin{figure}[htb]
\centering
\begin{tikzpicture}
    \cuboid{0}{2}{3}{0.5}{5}{5}{white}
    \cuboid{4}{0}{3}{5}{0.5}{5}{white}
    \cuboid{4}{2}{3}{5}{0.5}{5}{white}
    \cuboid{4}{4}{3}{5}{0.5}{5}{white}
    \path [every edge/.append style={draw=blue, |-|}]
        (1.5,-.5,5.5) edge ["$x_1$"'] (6.5,-.5,5.5)
        (6.6,-.35,5.5) edge ["$y_1$"'] (6.6,-.35,0.5)
        (6.65,-.25,0.5) edge ["$z_1$"'] (6.65,.25,0.5)
        (-.25,-1,5.5) edge ["$X$"']
        (6.5,-1,5.5);
    \node at (4,0,3){$M_1$};
    \node at (4,2,3){$M_2$};
    \node at (4,4,3){$M_3$};
    \node at (0,2,3){$M_4$};
    \end{tikzpicture}
    \caption{An arrangement of four modules. }\label{fig:floorplanning}
\end{figure}
For modules $i = 1,2,3$ orthogonal to the $Z$-axis, the thickness is $t_i = z_i$ and the area is $a_i = x_i y_i$. On the other hand, for module $4$ aligned orthogonal to the $X$-axis, we have $t_4 = x_4$ and $a_4 = y_4z_4$. Combined with the other constraints, the complete floorplanning problem is
\begin{equation}\label{eq:floorplanning}
\begin{aligned}
\text{minimize} \quad & \alpha XYZ + (1-\alpha)\sum_{i=1}^3 P_i K_i^{-1}z_i x_i^{-1}y_i^{-1} \\
&\qquad + (1-\alpha) P_4 K_4^{-1} x_4 y_4^{-1} z_4^{-1} \\
\text { subject to} \quad & x_i \geq x^{\min}_i, \quad y_i \geq y^{\min}_i,
\quad z_i \geq z^{\min}_i, \\ & i=1, \dots, 4,  \quad Z \leq Z^{\max}, \\
&\max\{x_1+x_4,x_2+x_4,x_3+x_4\} \leq X, \\
&\max\{y_1,y_2,y_3,y_4\} \leq Y,\\
&\max\{z_1+z_2+z_3,z_4\}\leq Z.
\end{aligned}
\end{equation}

We compare the computed results of \eqref{eq:floorplanning} for the 3D arrangement in Figure~\ref{fig:floorplanning} with a 2D arrangement when all four modules lie on the same plane. Using randomly generated parameters, the results are plotted in Figure~\ref{fig:floorplanning_experiment}. The 3D arrangement consistently outperforms the 2D arrangement, and the advantage increases as $\alpha$ increases, i.e., when the temperature is of lesser significance. 

\begin{figure}
  \includegraphics[trim={2em 5.8ex 1.9em 4.5ex},clip,width=.95\textwidth]{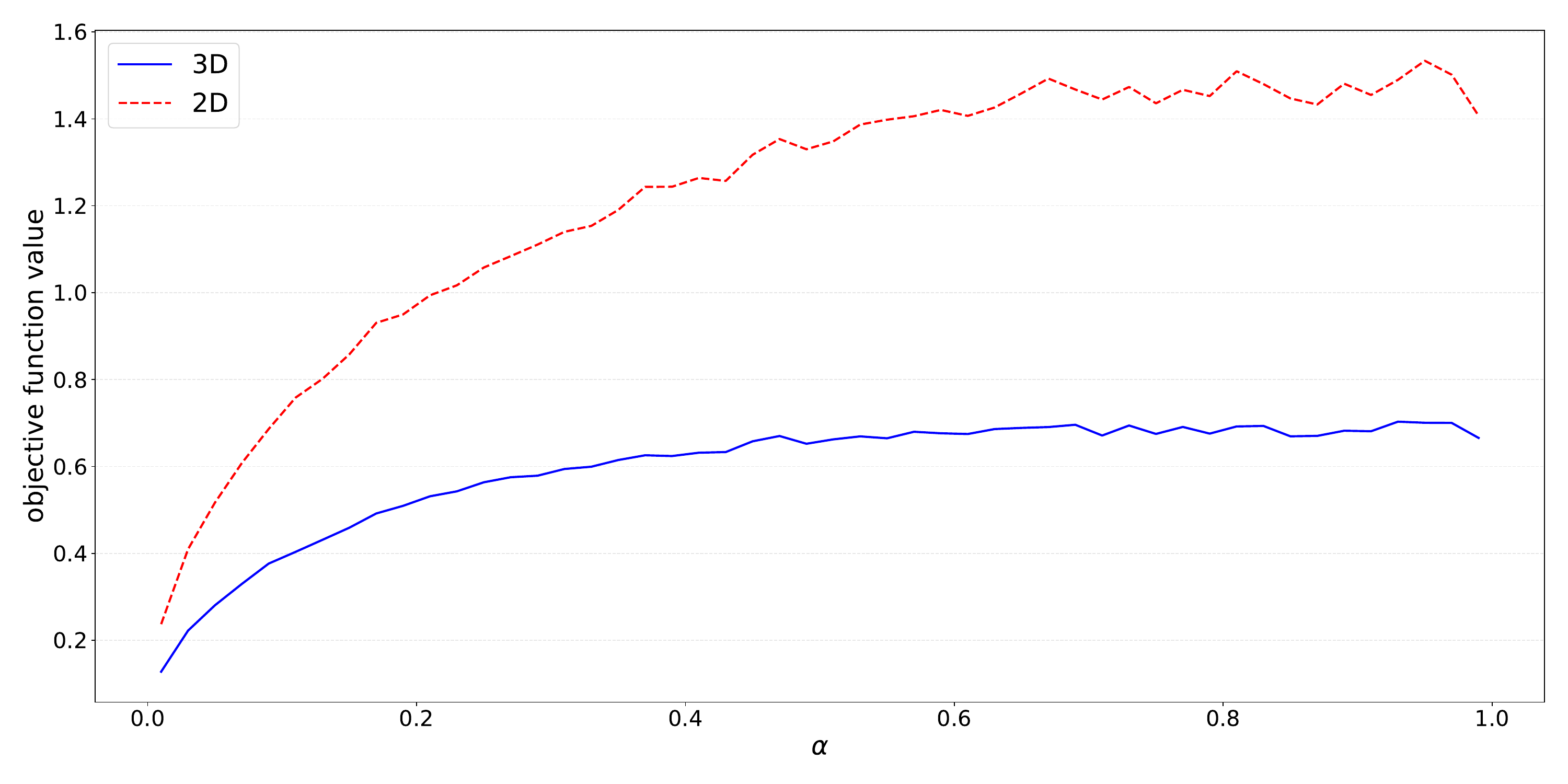}
  \caption{3D arrangement vs. 2D arrangement.}
  \label{fig:floorplanning_experiment}
\end{figure}
\end{example}

\begin{figure}
  \includegraphics[trim={2.1em 5.8ex 1.9em 4.5ex},clip,width=.95\textwidth]{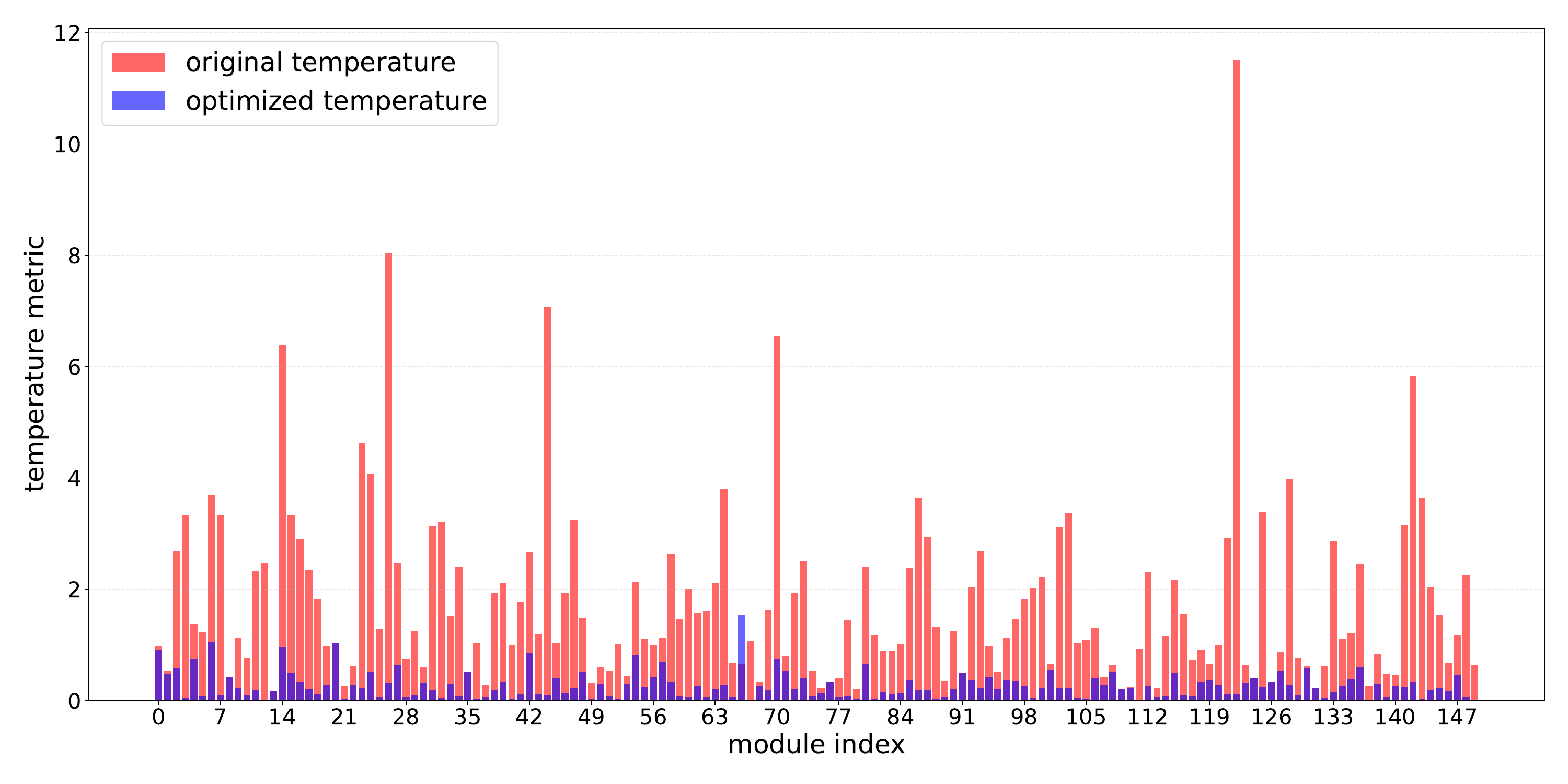}
  \caption{Temperature optimization of floorplanning.}
  \label{fig:floorplanning_temp}
\end{figure}

Moreover, the impact of transistor sizing on temperature becomes more pronounced as the number of modules increases. In Figure~\ref{fig:floorplanning_temp}, we conduct an experiment using a randomly generated 3D arrangement of 150 modules, with $\alpha = 0.6$. Due to the current lack of standard benchmarks for 3D modules, the experiment uses randomly generated parameters. Comparing module temperatures before and after optimization reveals a substantial temperature reduction across the modules, highlighting the effectiveness of the GP formulation in thermal management.

\subsection{Transistor sizing}\label{sec:sizing}
In electronic circuits, smaller transistors save space and reduce driver loads, while larger transistors can carry heavier loads to switch binary signals faster. Transistor sizing optimizes this trade-off to minimize delay. Gate sizing is a special case of transistor sizing when all transistors in a gate are scaled uniformly. We now introduce the model to formulate the full 3D IC transistor sizing problem analogous to the standard 2D approach \cite{digital_gp,transistor_sizing}.

Suppose the circuit topology $\mathcal{G} = (\mathcal{V}, \mathcal{E})$ is represented as a directed graph with transistors $1,\dots,n$. The following sets are determined by modeling the circuit as a resistor-capacitor (RC) model: 
\begin{enumerate}
    \item fan-in at $i$: $\operatorname{FI}(i) = \{j \in \mathcal{V}: (j,i) \in \mathcal{E}\}$,
    \item fan-out at $i$: $\operatorname{FO}(i) = \{j \in \mathcal{V}: (i,j) \in \mathcal{E}\}$,
    \item primary inputs: $\operatorname{PI} = \{i \in \mathcal{V} : \operatorname{FI}(i) = \varnothing\}$,
    \item primary outputs: $\operatorname{PO} = \{i \in \mathcal{V} : \operatorname{FO}(i) = \varnothing\}$,
    \item combinational logic block: $\operatorname{CB} = \mathcal{V} \setminus (\operatorname{PI} \cup \operatorname{PO})$.
\end{enumerate}
For each $i \in \operatorname{CB}$, we assign a scaling factor $x_i \geq 1$, which represents the optimization variables. When $x_i = 1$, the transistor is at its minimal size, referred to as unit scaling. A maximal scaling constraint $x_i \leq x_i^{\max}$ may also be added due to fabrication constraints.

Let the volume of gate $i$ at unit scaling be given by $\overline{\vol}_i$. The total volume of the IC is 
\[\vol = \sum_{i \in \operatorname{CB}} \overline{\vol}_i x_i,\]
a posynomial in $x$. 

We now use the resistor-capacitor (RC) delay model to formulate the GP. For each $i \in \operatorname{CB}$, the resistance of $i$ is
\[R_i = \overline{R}_i x^{-1},\]
where $\overline{R}_i$ is the resistance at unit scaling, given by the physical properties of the gate. 

If $i \in \operatorname{PO}$, then the input capacitance $C^{\operatorname{in}}_i$ is assumed given. For $i \in \operatorname{CB}$, the input capacitance and internal capacitance at unit scaling are given parameters, denoted $\overline{C}^{\operatorname{in}}_i$ and $\overline{C}^{\operatorname{int}}_i$. The input capacitance and internal capacitance are 
\[C^{\operatorname{in}}_i = \overline{C}^{\operatorname{in}}_i x_i \quad \text{and} \quad C^{\operatorname{int}}_i = \overline{C}^{\operatorname{int}}_i x_i.\]
The load capacitance $C^{\operatorname{L}}$ at $i \in \operatorname{CB} \cup \operatorname{PI}$ is
\[C^{\operatorname{L}}_i = \sum_{\mathclap{j \in \operatorname{FO}(i)}} C_j^{\operatorname{in}},\]
a posynomial in $x_i$, $i=1,\dots,n$. The total capacitance at $i \in \operatorname{CB} \cup \operatorname{PI}$ is given by
\[C_i = \begin{cases}
    C_i^L &\text{if } i \in \operatorname{PI},\\
    C_i^L + C_i^{\operatorname{int}} &\text{if } i \in \operatorname{CB}.
\end{cases}\]

The power of the circuit is  
\begin{equation}\label{eq:power}
  P = \sum_{i \in \operatorname{PI}\cup \operatorname{CB}} F_i C_i V_{dd}^2 + \sum_{i \in \operatorname{CB}} x_i \overline{I}_i V_{dd},
\end{equation}
where the unit scaling leakage current $\overline{I}_i$, activity frequency $F_i$, and supply voltage $V_{dd}$ are given.

A path in a circuit topology is a tuple $\operatorname{P} = (v_1,\dots, v_p)$ of circuit elements where $(v_i, v_{i+1}) \in \mathcal{E}$, $i = 0,\dots, p$, for some $v_0 \in \operatorname{PI}$ and $v_{m+1} \in \operatorname{PO}$. We denote the set of all paths in the circuit topology by $\mathcal{P}$. For $\operatorname{P} = (i_1,\dots,i_p) \in \mathcal{P}$, the delay at a gate $i_j \in \operatorname{P} \cap \operatorname{CB}$ and the delay of the path are respectively
\[
D_{i_j} = 0.69R_{i_j} C_{i_j} \quad \text{and} \quad D_{\operatorname{P}} =\sum_{j=2}^{p-1} D_{i_j}.
\]
The worst-case delay of the circuit is the maximal delay of all paths 
\[D = \max_{\operatorname{P} \in \mathcal{P}} D_{\operatorname{P}}. \]

A transistor sizing problem is a GPP of the form
\begin{equation}\label{eq:sizing}
  \begin{aligned}
      \text{minimize} \quad & D \\
  \text { subject to} \quad & P \leq P^{\max}, \quad \vol \leq \vol^{\max}, \\
  & 1 \leq x_i \leq x_i^{\max}, \quad i=1, \dots, n.
  \end{aligned}
\end{equation}

To illustrate the potential trade-offs between volume and delay, we conduct an illustrative experiment using the circuit topology in Figure~\ref{fig:nonplanar}. With randomly generated coefficients, we fix $P^{\max}$ and increase $\vol^{\max}$, plotting it against the corresponding optimal delay. As illustrated in Figure~\ref{fig:transistor_sizing}, the optimal delay follows an inverse trade-off against the maximal volume.

\begin{figure}[htb]
    \includegraphics[trim={1.7em 5ex 1.8em 4.5ex},clip,width=0.95\textwidth]{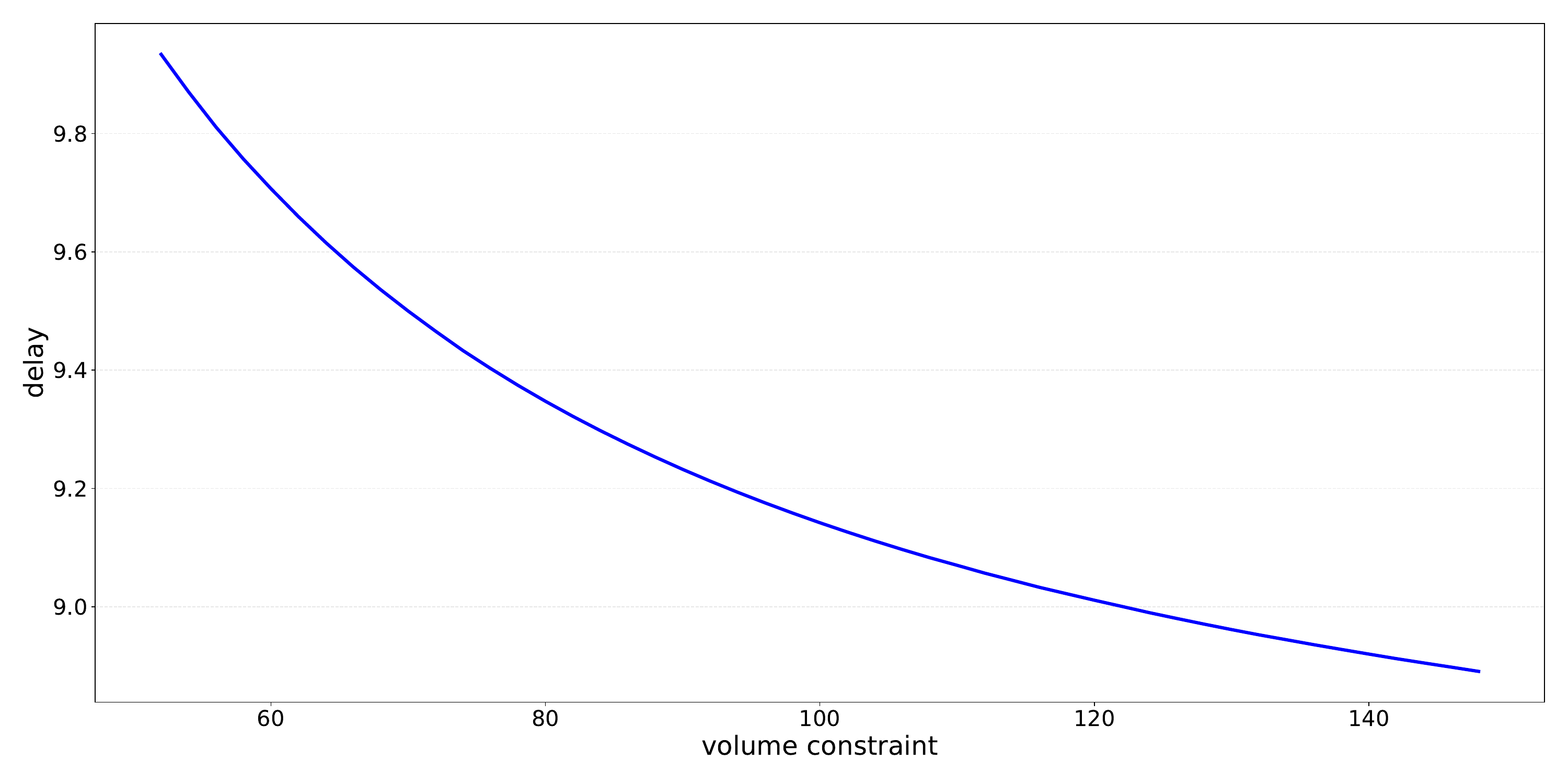}
    \caption{Volume vs delay trade-off. }\label{fig:transistor_sizing}
\end{figure}

To further illustrate the practical effectiveness of transistor sizing, we conduct gate-sizing experiments on circuits from the ISCAS-85 benchmark \cite{ISCAS85} using realistic parameters extracted from the ASAP7 predictive process design kit (PDK) for 7-nm FinFET technology \cite{ASAP7}.

\begin{figure}[htb]
    \includegraphics[trim={1.7em 4ex 1.8em 4.5ex},clip,width=0.95\textwidth]{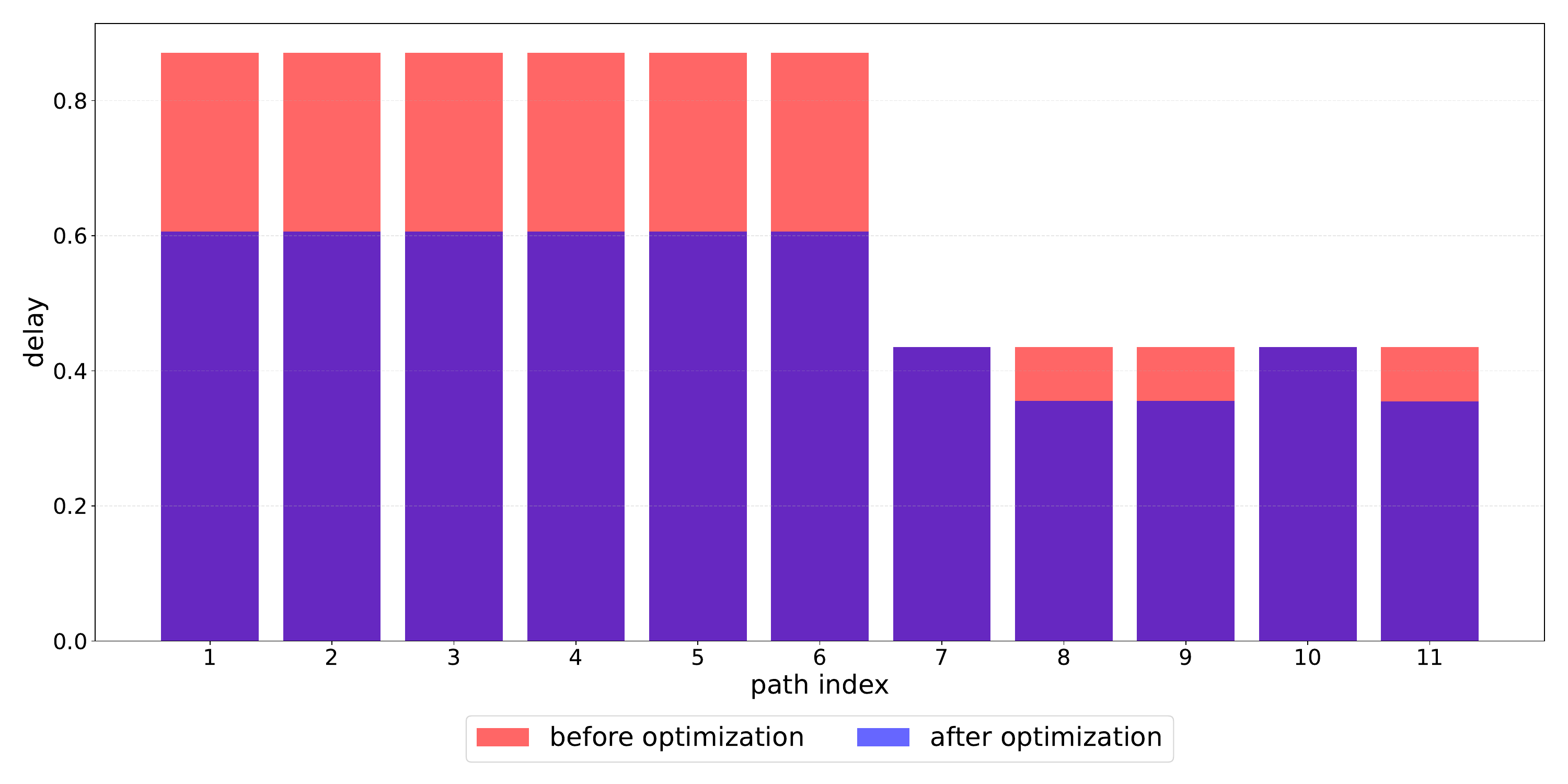}
    \caption{Gate sizing for c17 circuit.}\label{fig:c17}
\end{figure}

First, we consider the c17 circuit from ISCAS-85, which contains six NAND2 gates, each with four transistors, and $11$ distinct paths. In Figure~\ref{fig:c17}, we present a bar chart comparing the delays along each path before and after gate sizing optimization. The plot clearly demonstrates that delay is consistently reduced or maintained across all paths, with no degradation observed. Specifically, we achieve an overall improvement of approximately $30.36\%$. 

\begin{figure}[htb]
    \includegraphics[trim={1.7em 5ex 1.8em 4.5ex},clip,width=0.95\textwidth]{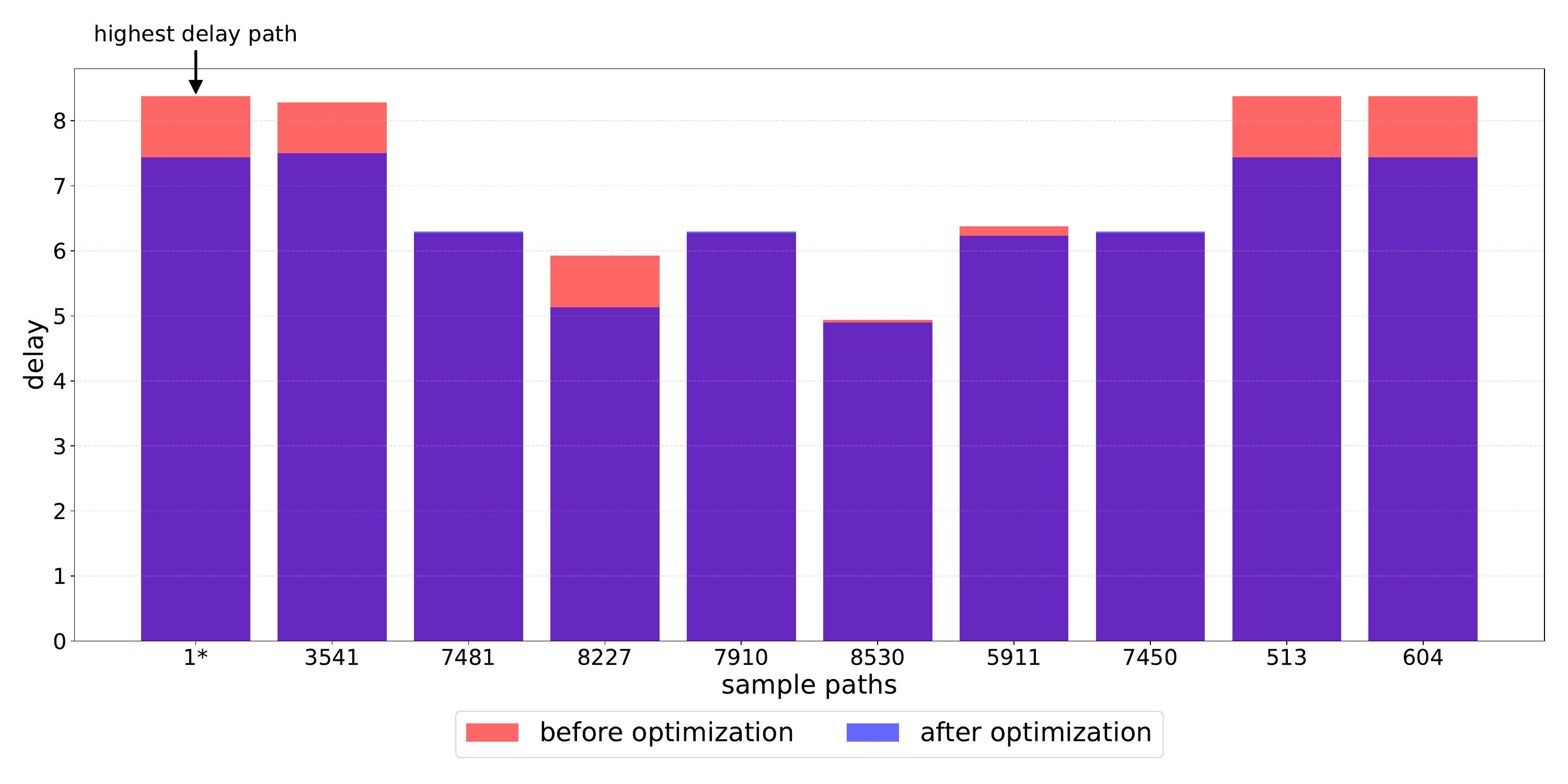}
    \caption{Gate sizing for c499 circuit. }\label{fig:c499}
\end{figure}

We apply the same methodology to a significantly larger benchmark circuit, c499 from ISCAS-85, which contains a total of $9440$ distinct paths. Due to the large number of paths, Figure~\ref{fig:c499} shows the delays of ten selected paths for clarity: one path that coincidentally has both the highest initial and final delay, and nine additional randomly sampled paths. The comparison reveals a reduction in delay across these representative paths, with the highest-delay path notably benefiting from optimization. The overall delay improvement for this circuit is approximately $11.24\%$. 

\begin{figure}[htb]
    \includegraphics[trim={1.7em 5ex 1.8em 4.5ex},clip,width=0.95\textwidth]{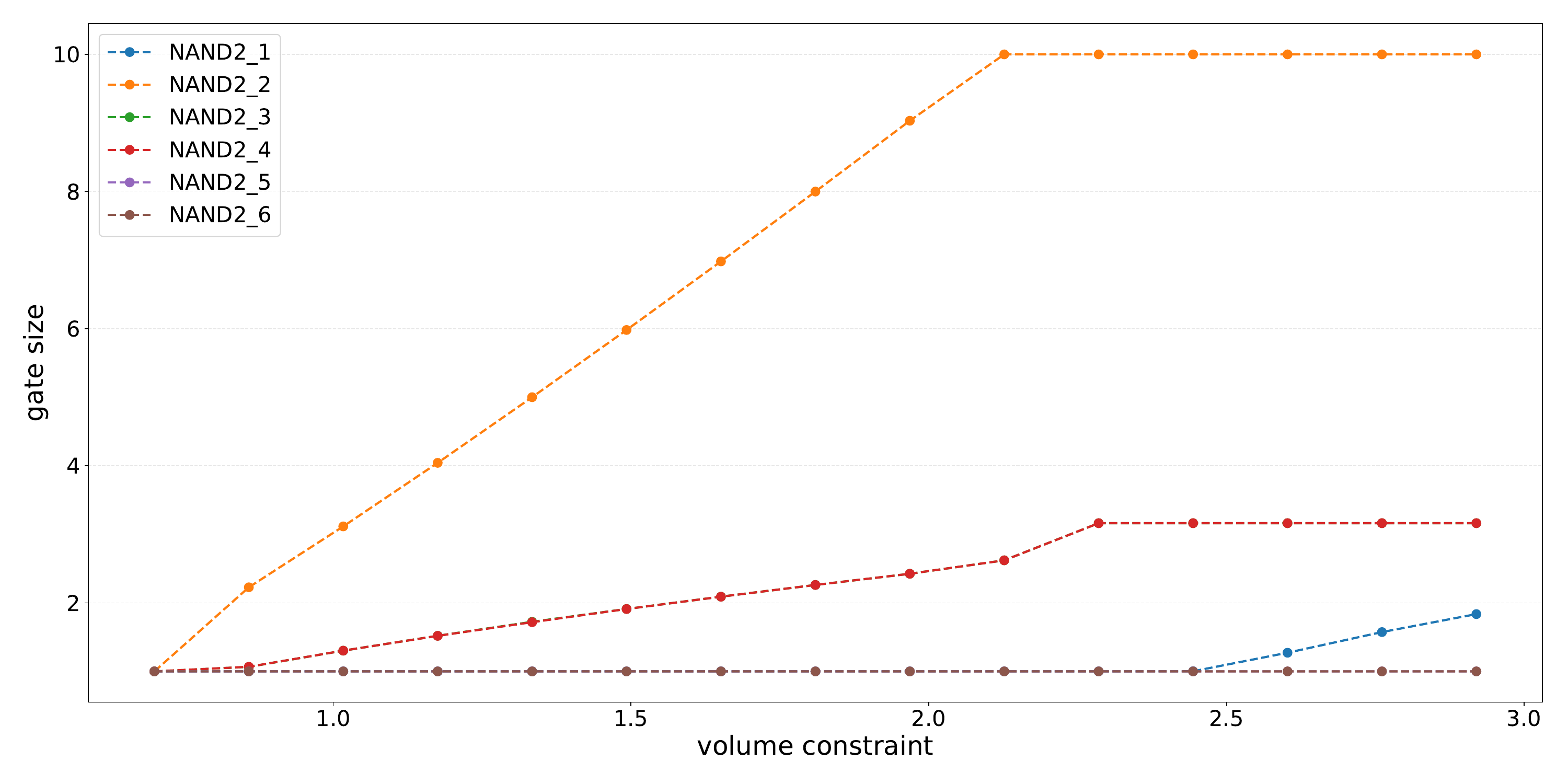}
    \caption{Individual gate sizes in c17 as volume constraint varies. }\label{fig:c17_individual}
\end{figure}

Next, we examine how individual gate sizes change as the volume constraint relaxes in c17. As shown in Figure~\ref{fig:c17_individual}, gate $2$ exhibits the most rapid growth in size compared to other gates, while gates $5$ and $6$, whose curves overlap in the plot, remain consistently small. This indicates that gate $2$ is critical to delay optimization, as its sizing has the greatest influence on the longest path delay. The divergence in gate sizes as the volume constraint increases highlights that transistor sizing optimization naturally prioritizes certain gates, selectively allocating resources to the most delay-sensitive parts of the circuit.

Formulation \eqref{eq:sizing} can be extended to include more detailed models, even for traditional MOSFETs, by incorporating additional constraints and factors \cite{digital_gp}. The complexity increases in full 3D designs with more detailed modeling of FinFETs. For example, external capacitance can be added, modeled as
\begin{equation}\label{eq:C_ext}
  C^{\operatorname{ext}} =\alpha \ln (1+ \alpha_2 x) + \alpha_3
\end{equation}
where $\alpha_1,\dots,\alpha_3$ are constants determined by physical properties of the FinFET \cite{fringing}. Although \eqref{eq:C_ext} is not strictly a generalized posynomial, it can be incorporated into GP extensions using appropriate approximations \cite{tutorial}. Given the dependency on specific circuit elements, a general approach to formulating these extensions of \eqref{eq:sizing} is unattainable.

\subsection{Interconnect sizing}\label{sec:sizing2}

Interconnect sizing generalizes wire sizing from planar ICs \cite{tutorial, wire_sizing} to full 3D ICs by incorporating both horizontal and vertical components. 3D Interconnects can be wires, vias, or other novel interconnect devices depending on the circuit type. Interconnect sizing determines the optimal length $l_i$ and width $w_i$ for each interconnect $i=1,\dots, m$ to minimize delay under fabrication constraints. 

\begin{figure}[htb]
    \resizebox{0.5\textwidth}{!}{\begin{circuitikz}
        \node[ground] at (0,0) {};
        \node[ground] at (3,0) {};
        \node[ground] at (6,0) {};
        \node[ground] at (6,3) {};
        \node[ground] at (9,6) {};
        \node[ground] at (9,3) {};
        \draw (0,0) to [american voltage source, invert, v=$V_{in}$](0,2) to [R = $R_1$] (0,4) to [generic = $1$] (3,4) -- (3,2) to [C = $C^{\operatorname{L}}_1$](3,0);
        \draw (3,2) to [generic = $5$] (6,2) to [C = $C^{\operatorname{L}}_5$](6,0);
        \draw (3,4) -- (3,6) to [generic = $2$] (6,6) -- (6,5) to [C = $C^{\operatorname{L}}_2$] (6,3);
        \draw (6,5) to [generic = $4$] (9,5) to [C = $C^{\operatorname{L}}_4$] (9,3);
        \draw (6,6) -- (6,8) to [generic = $3$] (9,8) to [C = $C^{\operatorname{L}}_3$] (9,6);
    \end{circuitikz}}
    \caption{An interconnect network.}\label{fig:network}
\end{figure}
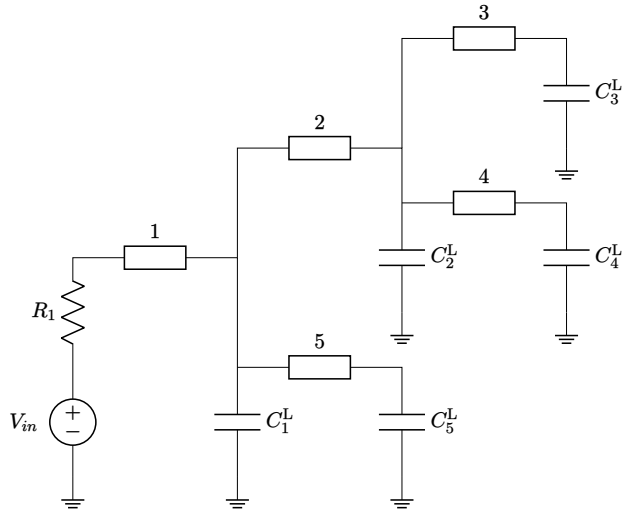

The interconnects are given in a tree-structured network where the input voltage serves as the root. Let $\mathcal{L}$ denote the set of leaves. For each interconnect $i$, let $\operatorname{P}(i)$ denote the unique path from the input voltage to $i$. For example, in Figure~\ref{fig:network}, $\mathcal{L} = \{3, 4, 5\}$, $\operatorname{P}(5) = (1, 5)$, and $\operatorname{P}(4) = (1, 2, 4)$. 

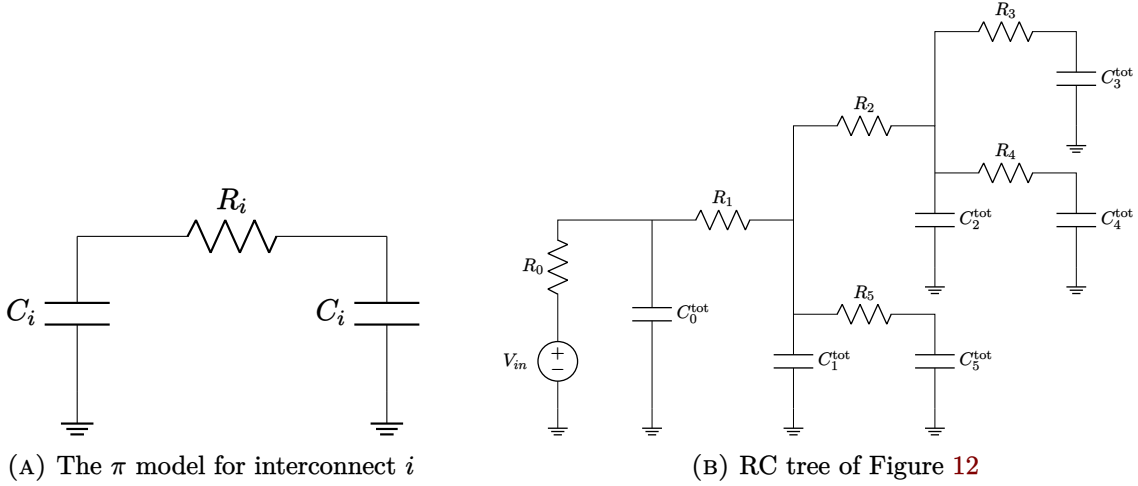
\begin{figure}
    \begin{subfigure}{0.45\textwidth}
      \centering
      \resizebox{0.75\textwidth}{!}{%
        \begin{circuitikz}
          \node[ground] at (-2,0) {};
          \node[ground] at (2,0) {};
          \draw (-2,0) to [C = $C_i$](-2,2);
          \draw (2,0) to [C = $C_i$](2,2);
          \draw (-2,2) to [R = $R_i$](2,2);
        \end{circuitikz}}
      \caption{The $\pi$ model for interconnect $i$}
      \label{fig:pi}
    \end{subfigure}
    \begin{subfigure}{0.54\textwidth}
      \resizebox{0.95\textwidth}{!}{%
        \begin{circuitikz}
          \node[ground] at (-2,0) {};
          \node[ground] at (0,0) {};
          \node[ground] at (3,0) {};
          \node[ground] at (6,0) {};
          \node[ground] at (6,3) {};
          \node[ground] at (9,6) {};
          \node[ground] at (9,3) {};
          \draw (0,4) to [C = $C^{\operatorname{tot}}_0$](0,0);
          \draw (-2,0) to [american voltage source, invert, v=$V_{in}$](-2,2) to [R = $R_0$] (-2,4) to (0,4) to [R = $R_1$] (3,4) -- (3,2) to [C = $C^{\operatorname{tot}}_1$](3,0);
          \draw (3,2) to [R = $R_5$] (6,2) to [C = $C^{\operatorname{tot}}_5$](6,0);
          \draw (3,4) -- (3,6) to [R = $R_2$] (6,6) -- (6,5) to [C = $C^{\operatorname{tot}}_2$] (6,3);
          \draw (6,5) to [R = $R_4$] (9,5) to [C = $C^{\operatorname{tot}}_4$] (9,3);
          \draw (6,6) -- (6,8) to [R = $R_3$] (9,8) to [C = $C^{\operatorname{tot}}_3$] (9,6);
        \end{circuitikz}}
      \caption{RC tree of Figure~\ref{fig:network}}
      \label{fig:RCtree}
    \end{subfigure}
    \caption{A $\pi$ model and the RC tree of Figure~\ref{fig:network}.}\label{fig:interconnect}
\end{figure}

To compute the delay, we model the interconnect network as an RC tree using a $\pi$ model for each interconnect \cite{tutorial}. The $\pi$ model replaces each interconnect by a $\pi$-shaped circuit with two capacitors and one resistor as illustrated in Figure~\ref{fig:interconnect}. The resistance and capacitance in the $\pi$ model are respectively
\[R_i = \alpha_i l_i w_i^{-1}, \quad C_i = \beta_i l_i w_i + \gamma_i l_i, \quad i = 1,\dots, m,\]
where $\alpha_i$, $\beta_i$, and $\gamma_i$ are given positive constants determined by the physical properties of each interconnect. The function $R_i$ is a monomial, while $C_i$ is a polynomial in $l_i$ and $w_i$ for each $i=1, \dots, m$. 

In the RC tree, the downstream elements $\operatorname{DS}(i)$ are those appearing to the right of $i$. Assuming the capacitive load $C^{\operatorname{L}}_i$ for each interconnect is given, the total capacitance $C^{\operatorname{tot}}_i$ is the sum of $C^{\operatorname{L}}_i$, $C_i$, and $C_j$ for all elements immediately downstream from $i$. For instance, in Figure~\ref{fig:RCtree}, $C^{\operatorname{tot}}_1 = C_1^{\operatorname{L}} + C_1 + C_2 + C_5$ and $C^{\operatorname{tot}}_2 = C_2^{\operatorname{L}} + C_2 + C_3 + C_4$. Since $C^{\operatorname{tot}}_i$ are sums of $C_i$ and positive constants, they are posynomials in $l_i$ and $w_i$, $i=1, \dots, m$.

We use the Elmore delay model \cite{VLSI_book} to formulate the GP. The total delay from the root to capacitor at interconnect $i$ is 
\[D_i = \sum_{j \in \operatorname{P}(i)} R_j \Biggl(\sum_{k \in \operatorname{DS}(j)} C^{\operatorname{tot}}_k\Biggr),\]
a posynomial for each $i=1,\dots,m$. The interconnect sizing problem is thus:
\begin{equation*}
  \begin{aligned}
      \text{minimize} \quad & \max_{i \in \mathcal{L}} D_i \\
  \text { subject to} \quad & w_i^{\min} \leq w_i \leq w_i^{\max}, \quad l_i^{\min} \leq l_i \leq l_i^{\max}, \\
  & i=1, \dots, m, \quad \sum_{i=1}^m l_iw _i^2 \leq \vol^{\max}.
  \end{aligned}
\end{equation*}
where $w_i^{\min}$, $w_i^{\max}$, $l_i^{\min}$, $l_i^{\max}$ are fabrication constraints on the interconnects and $\vol^{\max}$ is a predetermined cap on the total volume.

\begin{figure}
    \centering
    \begin{subfigure}{0.49\textwidth}
    \centering
    \includegraphics[trim={2em 1.6ex 1.9em 4.4ex},clip,width=\textwidth]{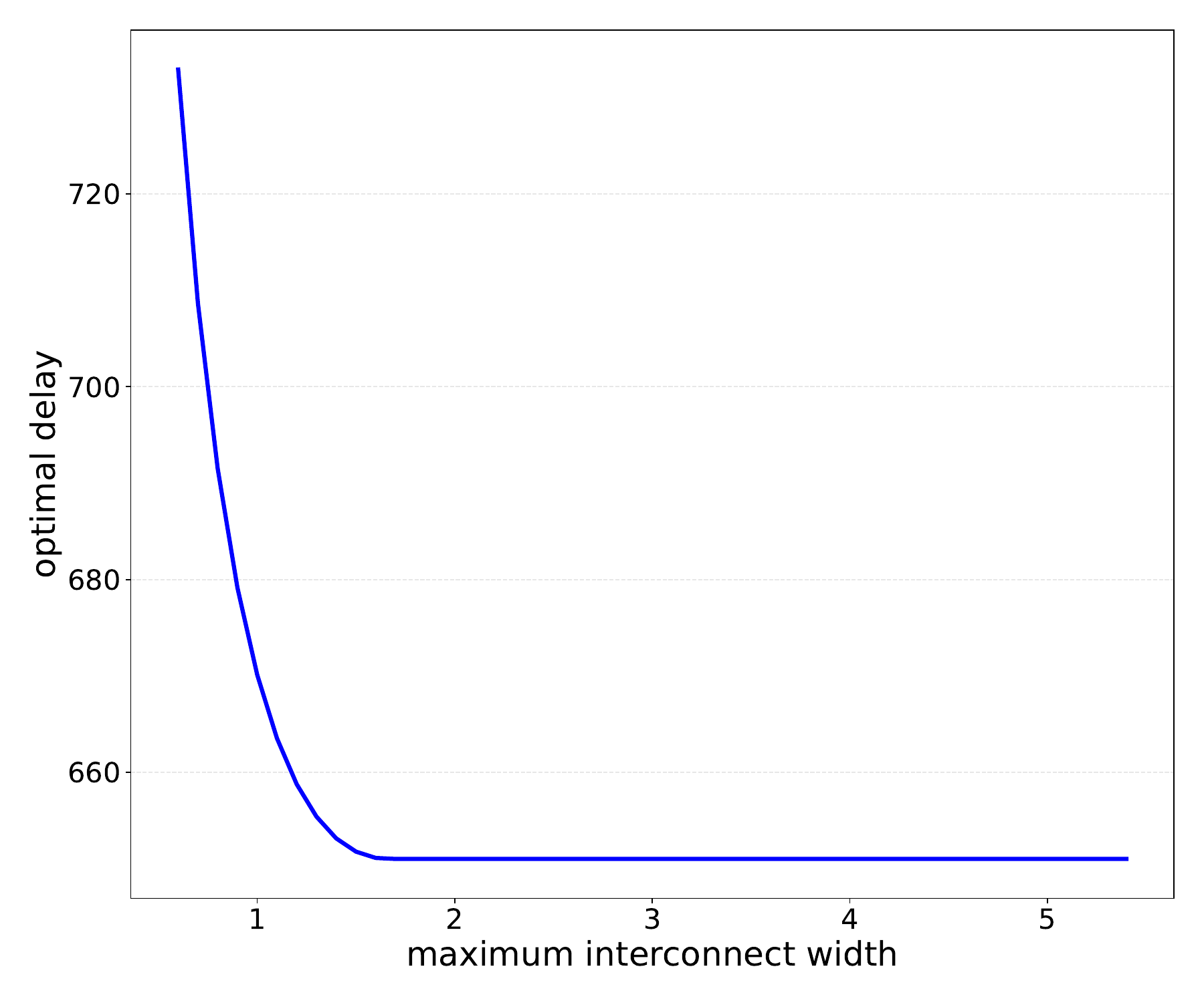}
    \caption{Small interconnect width}
    \label{fig:transistor_sizing2}
      \end{subfigure}
    \begin{subfigure}{0.49\textwidth}
        \centering
      \includegraphics[trim={2em 1.6ex 1.9em 4.4ex},clip,width=\textwidth]{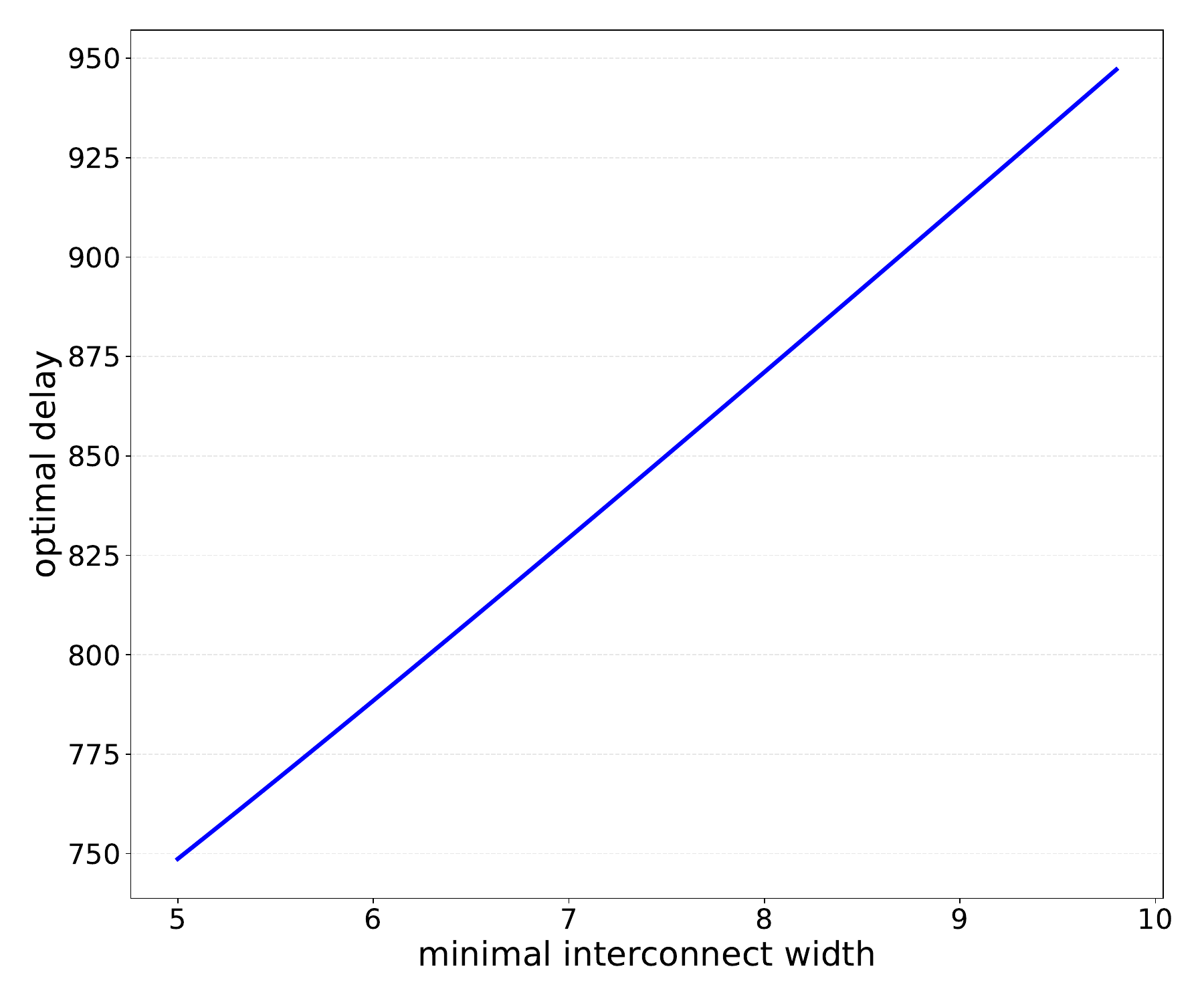}
    \caption{Large interconnect width}
    \label{fig:transistor_sizing}
    \end{subfigure}
    \caption{Dependence of delay on interconnect width.}\label{fig:interconnect_width}
\end{figure}

We investigate how interconnect width influences optimal delay using the interconnect network depicted in Figure~\ref{fig:network}. In our first experiment, we assume all interconnects have similar dimensions and randomly generate the coefficients. By incrementally increasing the upper bound $w_1^{\max}$, we plot the corresponding optimal delay, as shown in Figure~\ref{fig:transistor_sizing2}. The optimal delay initially decreases as $w_1^{\max}$ grows, but eventually stabilizes at approximately $1.66109708$. This stabilization occurs due to the behavior of the optimal interconnect width $w_1^*$, which, according to our computations, satisfies:
\[
w_1^* = \begin{cases}
    w_1^{\max} & \text{if } w_1^{\max} \leq 1.66109708, \\
    1.66109708 & \text{if } w_1^{\max} > 1.66109708.
\end{cases}
\]
This experiment empirically identifies $1.66109708$ as a critical threshold for interconnect width, beyond which further increases do not yield additional improvements in delay optimization.

Next, we consider the scenario where one interconnect (for example, a TSV) is significantly larger than the others. With all other coefficients fixed, we systematically vary the lower bound $w_1^{\min}$ and observe its effect on the optimal delay. In this scenario, the optimal width consistently occurs at $w_1^* = w_1^{\min}$. Furthermore, as illustrated in Figure~\ref{fig:transistor_sizing}, the optimal delay increases linearly with $w_1^{\min}$.

\begin{figure}
    \centering
    \begin{subfigure}{0.49\textwidth}
    \centering
    \includegraphics[trim={2em 1.6ex 1.9em 4.4ex},clip,width=\textwidth]{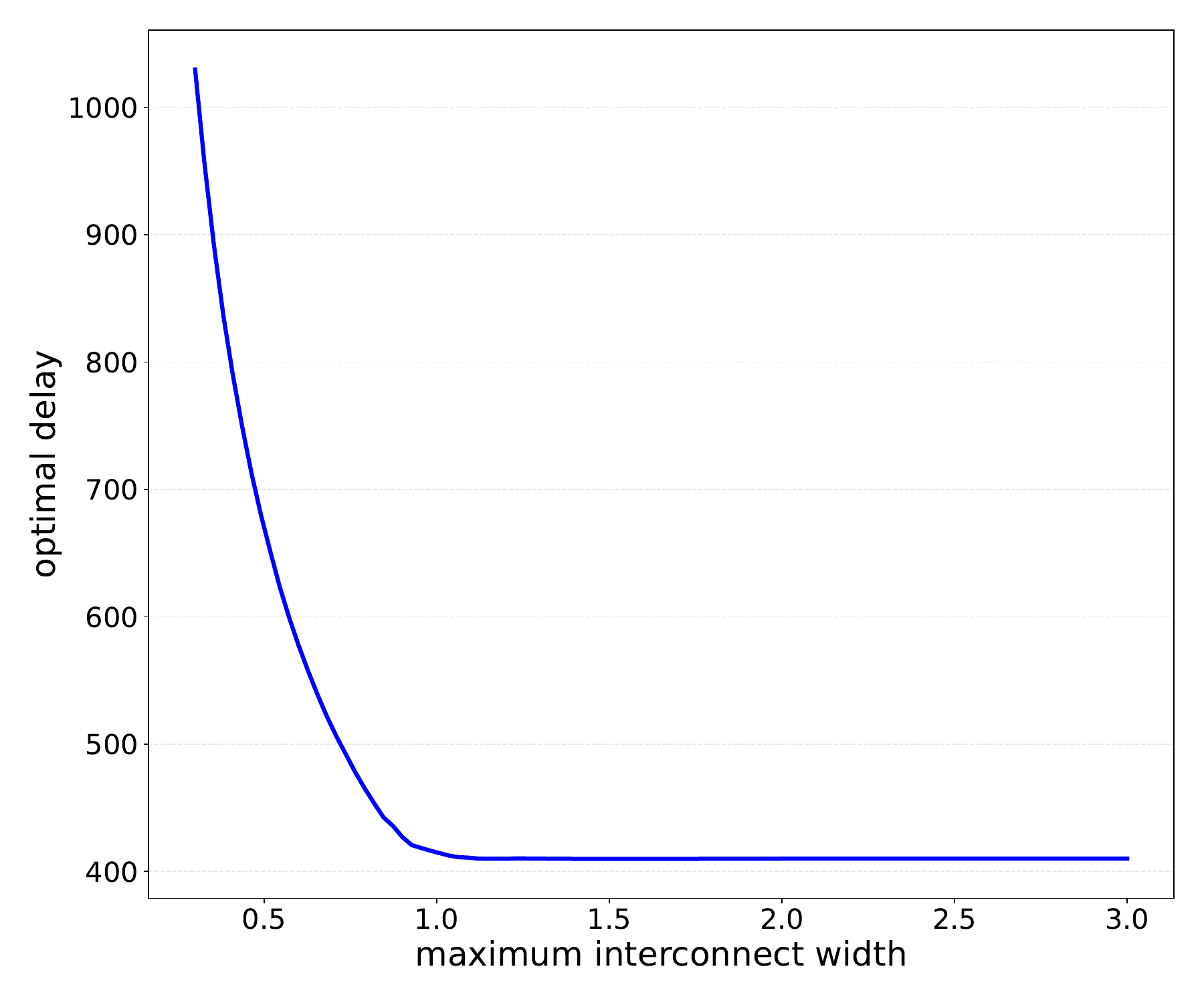}
    \caption{Decreasing upper bound}
    \label{fig:c17max}
      \end{subfigure}
    \begin{subfigure}{0.49\textwidth}
        \centering
      \includegraphics[trim={2em 1.6ex 1.9em 4.4ex},clip,width=\textwidth]{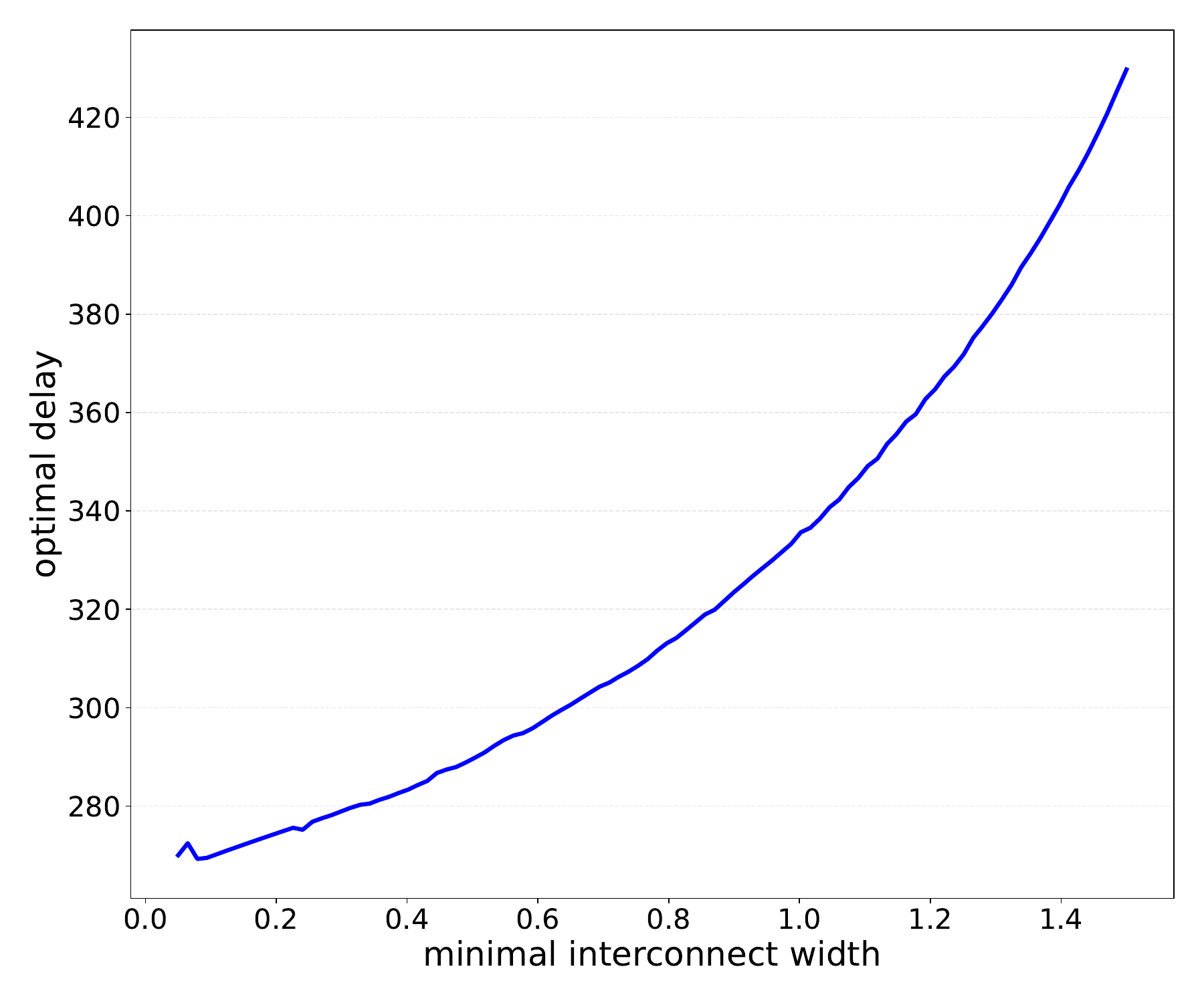}
    \caption{Increasing lower bound}
    \label{fig:c17min}
    \end{subfigure}
    \caption{Dependence of delay on interconnect width for c17 circuit.}\label{fig:c17_width}
\end{figure}

This observed trend is general rather than coincidental. To validate this behavior, we conduct additional experiments using the interconnect network derived from the \texttt{c17} circuit of the ISCAS-85 benchmark \cite{ISCAS85}. To highlight the effect, we use a uniform upper bound across all interconnects:
\[w^{\max} = w_1^{\max} = \cdots = w_n^{\max}\]
and randomly generate other parameters. Similarly, in a separate experiment, we apply a uniform lower bound:
\[w^{\min} = w_1^{\min} = \cdots = w_n^{\min}\]
again randomly generating the remaining parameters. Figure~\ref{fig:c17max} confirms the stabilization phenomenon observed previously in Figure~\ref{fig:transistor_sizing2}. However, since the lower bounds simultaneously increase for all interconnects, the delay now increases polynomially rather than linearly, reflecting the compounded effects across multiple interconnects.

\section{Conclusion}

This article is likely the first systematic work to consider mathematical optimization problems in full 3D circuit integration. The graph embedding algorithms and GP-based optimization methods introduced in this paper are efficient, systematic, and accessible. The primary limitation of this work is that our formulations are based on current technological understanding. As new engineering advancements emerge, future work will need to address the challenges posed by these innovations. Such advancements may introduce novel optimization objectives and modeling terms for GP formulations or shift design priorities beyond minimizing bends in the layout problem.

\bibliographystyle{abbrv}
\bibliography{3dcirc}
\end{document}